\documentclass[12pt,amsfonts]{article}
\usepackage[margin=1in]{geometry}  
\usepackage{graphicx}              
\usepackage{amsmath}               
\usepackage{amsfonts}              
\usepackage{amsthm}                
\def\nd{\noindent}
\def\thend{\rule{3mm}{3mm}}

\newtheorem{claim}{Claim}[section]
\newtheorem{thm}{Theorem}[section]

\newtheorem{lem}{Lemma}[section]

\numberwithin{equation}{section}

\begin{document}

\title{Existence of heteroclinic solution for a class of non-autonomous second-order  equation.  }
\author{\sf Claudianor O. Alves\thanks{Research of C. O. Alves partially supported by  CNPq 304036/2013-7  and INCT-MAT} 
\\
\small{Universidade Federal de Campina Grande, }\\
\small{Unidade Acad\^emica de Matem\'atica } \\
\small{CEP:58429-900, Campina Grande-PB, Brazil}\\}
\date{}
\maketitle
\begin{abstract}
In this paper, we use variational methods to prove the existence of heteroclinic solutions for a class of non-autonomous second-order equation.
\end{abstract}

\vspace{0.5 cm}


\section{Introduction}

Consider the non-autonomous second-order differential equation
\begin{equation} \label{E1}
\ddot{x}(t)=a(\epsilon  t) V'(x(t)), \,\,\, t \in \mathbb{R},
\end{equation}
\begin{equation} \label{E2}
x(t) \to -1 \,\,\, \mbox{as} \,\,\, t \to -\infty, \,\, x(t) \to 1 \,\,\, \mbox{as} \,\,\, t \to +\infty,
\end{equation}
where $\epsilon >0$ is a positive parameter and $V:\mathbb{R} \to \mathbb{R}$ is a function verifying: \\

\noindent $ (V_1) \,\, V \in C^{2}(\mathbb{R}^{N}, \mathbb{R}).$ \\

\noindent $ (V_2) \,\, V(x) \geq 0 \,\, \forall x \in \mathbb{R}$ and $V(-1)=V(1)=0 .$ \\

\noindent $ (V_3) \,\,\,   V(x)>0$ \,\, for all $ x \in (-1,1),$ \\

\noindent and \\

\noindent $(V_4) \,\,\, V''(-1),V''(1)>0$. \\

Related to function $a: \mathbb{R} \to \mathbb{R}$, we assume that it is a bounded continuous function satisfying some conditions which will be mentioned later on. \\

The main goal of the present paper is to prove the existence of solution for problem (\ref{E1})-(\ref{E2}), which is called a {\it heteroclinic solution}, connecting the equilibria -1 and 1. \\

The existence of heteroclinic solution has received a special attention, because this type of solution appears in a lot of mathematical models, such as Mechanics, Chemistry and Biology, for more details about this subject, we cite Bonheure and Sanchez \cite{BS}.

In \cite{BS}, the existence of heteroclinic solution for (\ref{E1})-(\ref{E2}) has been studied for some classes of function $a$. More precisely, in that paper the following classes were considered:

\vspace{0.5 cm}

\noindent {\bf Class 1:} \,  $a$ is a positive constant. \\

\noindent {\bf Class 2:} \, $a$ is a periodic continuous function with 
$$
\inf_{t \in \mathbb{R}}a(t)=a_0>0. \eqno{(a_0)}
$$

\noindent {\bf Class 3:} \, $a$ is a bounded continuous function and there are $a_1,a_2>0$ verifying
$$
a_1 \leq a(t) \leq a_2 \,\,\, \forall t \in \mathbb{R}  \eqno{(a_1)}
$$  
and
$$
a(t) \to a_2, \,\,\, \mbox{as} \,\,\,\, |t| \to +\infty, \eqno{(a_2)}
$$
with $a(t)<a_2$ in some set of nonzero measure. 

\vspace{0.5 cm}

In \cite{GaSa}, Gavioli and Sanchez have assumed that $a$ belongs to ensuing class : \\ 

\noindent {\bf Classe 4:} \, There is $t_0$ such that $a$ is increasing in $(-\infty,t_0]$,  $a$ is decreasing in $[t_0,+\infty)$, \linebreak $\displaystyle \lim_{|t| \to +\infty}a(t)=l>0$, and
$$
\lim_{|t|\to +\infty}|t|(l-a(t))=0. \eqno{(a_3)}
$$

Gavioli in \cite{Ga1} has studied the following class \\

\noindent {\bf Class 5:} \, There are $0<l<L$ such that  
$$
l \leq a(t) \leq L \,\,\, \forall t \in \mathbb{R}, \eqno{(a_4)}
$$ 
$$
a(t) \to L \,\,\, \mbox{as} \,\,\, |t| \to +\infty, \eqno{(a_5)}
$$
and $L/l$ is suitably bounded from above. 

\vspace{0.5 cm}

After, Gavioli in \cite{Ga2} considered the situation where $a$ is in the class  \\

\noindent {\bf Class 6:} \, $a \in L^{\infty}(\mathbb{R},[0,+\infty))$ and there are $l>0, S<T$, such that 
$$
a(t)=l \,\,\, \mbox{for} \,\,\,\, t \notin [S,T]. \eqno{(a_6)}
$$ 

Finally, in \cite{S},  Spradlin established the existence of heteroclinic for the case where $a$ within class

\vspace{0.5 cm}

\noindent {\bf Class 7:} There are $\underline{l},l>$ such that
$$
a(t) \to l \,\,\, \mbox{as} \,\,\, |t| \to +\infty,
$$
and
$$
\underline{l} \leq a(t) \leq L= \nu\sqrt{l \underline{l}}/ \int_{-1}^{1}\sqrt{V(x)}\,dx,
$$
where
$$
\nu=\min\left\{\int_{-1}^{\xi_{-}}\sqrt{V(x)}\,dx, \int_{\xi_{+}}^{1}\sqrt{V(x)}\,dx  \right\}
$$
with
$$
\xi_{-}=\min \left\{ x:x>-1, \, V'(x)=0  \right\} \,\,\,\, \mbox{and} \,\,\,\, \xi_{+}=\max \left\{ x:x<1, \, V'(x)=0  \right\}.
$$

In all above references, the main idea to get a solution for (\ref{E1})-(\ref{E2}) is looking for critical point for the functional $J:H^{1}_{loc}(\mathbb{R}) \to [0,+\infty]$ given by
$$
J(x)=\int_{-\infty}^{+\infty}\left( \frac{1}{2}|\dot{x}|^{2}+a(\epsilon t){V}(x(t))  \right)dt.
$$
In some of the above references, the existence of critical point was established showing that $J$ possesses a critical point on one of the ensuing sets
$$
\Sigma=\left\{x \in H^{1}_{loc}(\mathbb{R}): x(-\infty)=-1 \,\,\, \mbox{and} \,\,\, x(+\infty)=1 \right\}
$$ 
or
$$
W=\left\{ x \in H^{1}_{loc}(\mathbb{R}): x+1 \in H^{1}((-\infty, 0]), \, x-1 \in H^{1}([0,+\infty))  \right\}.
$$
The main tool used is the variational method, more precisely, deformation lemma and minimization techniques.

Motivated by cited references, we intend to study the existence of heteroclinic solution for (\ref{E1})-(\ref{E2}) for three new classes of function $a$. Here, we will consider the  following classes: \\

\noindent {\bf Class 8:}\, $a \in L^{\infty}(\mathbb{R})$ and 
$$
\liminf_{|t| \to \infty}a(t)=a_\infty > \inf_{t \in \mathbb{R}}a(t)=a(0)>0. \eqno{(a_7)}
$$  
This class of functions was introduced by Rabinowitz \cite{R12} to study existence of solution for a P.D.E. of the type 
$$
-\Delta u +V(\epsilon x)u=f(u), \,\,\,\, \mathbb{R}^{N}.
$$
This way, throughout this article,  we will called it of Rabinowitz's condition.

\vspace{0.5 cm}

\noindent {\bf Class 9:}\, $a$ is asymptotically periodic, that is, there is a continuous periodic function $a_P:\mathbb{R} \to \mathbb{R}$ satisfying: 
$$
|a(t)-a_P(t)| \to 0 \,\,\, \mbox{as} \,\,\,\ |t| \to +\infty \eqno{(a_8)}
$$ 
and
$$
0<\inf_{t \in \mathbb{R}}a(t)\leq a(t)<a_P(t) \,\,\, \forall t \in \mathbb{R}.  \eqno{(a_9)}
$$

\vspace{0.5 cm}

\noindent {\bf Class 10:}\, $a$ is coercive, that is, 
$$
0<\inf_{t \in \mathbb{R}}a(t) \,\,\, \mbox{and} \,\,\,\, a(t) \to +\infty \,\,\,\, \mbox{as} \,\,\,\, |t| \to +\infty. \eqno{(a_{10})}
$$

Our main result is the following 

\begin{thm} \label{T1} Assume $(V_1)-(V_4)$ and that $a$ belongs to Class 9 or 10. Then, for each $\epsilon >0$, problem (\ref{E1})-(\ref{E2}) has a solution $x \in H_{loc}^{1}(\mathbb{R}) \cap C^{2}(\mathbb{R})$ and $x(t) \in (0,1)$ for all $t \in \mathbb{R}$. If $a$ belongs to  Class 8, the existence of solution is established for $\epsilon$ small enough. 

\end{thm}

In the proof of Theorem \ref{T1}, we explored some arguments used in \cite{BS} and \cite{S}. The basic idea is working with a minimization problem, which will lead us to get a heteroclinic solution for the problem (\ref{E1})-(\ref{E2}), for more details, see Sections 3, 4 and 5.

Before to conclude this introduction, we would like to cite the papers of Bonheure, Sanchez and Tarallo \cite{BSTT}, Bonheure, Obersnel and  Omari \cite{BOO}, Bonheure, Coelho and Nys \cite{BCN}, Coti Zelati and Rabinowitz \cite{CR}, Korman, Lazer and Li \cite{KLL}, Rabinowitz \cite{R2}, and their references, where the reader can find interesting results about the existence of heteroclinic solutions for related problems.

The plan of the paper is as follows: In Section 2, we prove some technical results, which will be useful in the proof of Theorem \ref{T1}. In Sections 3 and 4, we study the case where $a$ verifies the Rabinowitz's condition and it is asymptotically periodic respectively, while the coercive case is considered in Section 5. In Section 6, we make some final considerations.

\section{Technical results }

In this section, we will show some results, which are crucial in the proof of Theorem \ref{T1}. However, we would like to point out that in their proofs it is enough to assume that function $a$ verifies the following condition:

\vspace{0.5 cm}

There are $l_0, l_1 >0$ such that
$$
l_0 \leq a(t) \leq l_1 \,\,\,\, \forall t \in \mathbb{R}. \eqno{(a_{11})}
$$

To begin with, we observe that from $(V_1)-(V_4)$, there are $C_1,C_2, \delta >0$ with $C_1<C_2$, such that 
\begin{equation} \label{EST1}
C_1(x-1)^{2} \leq V(x) \leq C_2 (x-1)^{2} \,\,\,\,\,\,\, \forall x \in (1-\delta, 1+\delta)
\end{equation}
and
\begin{equation} \label{EST2}
C_1(x+1)^{2} \leq V(x) \leq C_2 (x+1)^{2} \,\,\,\,\,\,\, \forall x \in (-1-\delta, -1+\delta).
\end{equation}

In what follows, we will make a modification on function $V$, by assuming that it satisfies the following properties: \\

\noindent $ (V_5) \,\,\, V(x) >0$ \, for all  $ x \in (-\infty, -1-\delta) \cup (1+\delta, +\infty)$,\\

\noindent $ (V_6) \,\, V'(x)x > 0$ \,  for all  $ x \in (-\infty, -1) \cup (1, +\infty)$,  \\

\noindent and \\

\noindent $(V_7) \,\,\, V(x) \to +\infty $ \,\, as \,\,\, $|x| \to +\infty$. \\

Hereafter, we will denote by $\tilde{V}$ the new function. This way, 
$$
\tilde{V}(x)=V(x) \,\,\, \forall x \in (-1-\delta,1+\delta), \,\,\, \tilde{V}'(x)x >0 \,\, \mbox{for} \,\, |x| >1 \,\, \mbox{and} \,\, \tilde{V}(x) \to +\infty \,\,\, \mbox{as} \,\,\, |x| \to +\infty.
$$
Moreover, we denote by $W \subset H^{1}_{loc}(\mathbb{R})$ the set  
$$
W=\left\{ x \in H^{1}_{loc}(\mathbb{R}): x+1 \in H^{1}((-\infty, 0]), \, x-1 \in H^{1}([0,+\infty))  \right\},
$$
and by $J_\epsilon:H^{1}_{loc}(\mathbb{R}) \to [0,+\infty]$ the functional given by 
$$
J_\epsilon(x)=\int_{-\infty}^{+\infty}\left( \frac{1}{2}|\dot{x}|^{2}+a(\epsilon t)\tilde{V}(x(t))  \right)dt.
$$
From $(V_1)-(V_4)$, more precisely (\ref{EST1})-(\ref{EST2}), it follows that $J_\epsilon(x) < +\infty$ for all $x \in W$ and $J_\epsilon$ is Fr\'echet differentiable, in the sense that, 
$$
J_\epsilon'(x)v = \int_{-\infty}^{+\infty}\left( \dot{x}\dot{v}+a(\epsilon t)\tilde{V}(x(t))v(t) \right)dt \,\,\,\,\,\,\, \forall \, x \in W \,\,\, \mbox{and} \,\,\, v \in H^{1}(\mathbb{R}).
$$

In the sequel, we say that $(x_n)$ is a $(PS)_c$ sequence for $J_\epsilon$, if $(x_n) \subset W$ with  
$$
J_\epsilon(x_n) \to c \,\,\, \mbox{and} \,\,\ \|J_\epsilon'(x_n)\|_* \to 0 \,\,\, \mbox{as} \,\,\, n \to \infty,
$$
where
$$
\|J_\epsilon'(x)\|_*=\sup \left\{J_\epsilon'(x)v: \, v \in H^{1}(\mathbb{R}), \|v\| =1   \right\},
$$
and $\|\,\,\,\|$ denotes the usual norm in $H^{1}(\mathbb{R})$.

\vspace{0.5 cm}

The next two lemmas can be found in \cite{S}, however for reader's convenience we will write their proofs.  

\begin{lem} \label{L1} If $x \in H^{1}_{loc}(\mathbb{R})$ and $J_\epsilon(x)<\infty$, then
$$
x(t) \to -1 \,\,\, \mbox{or} \,\,\, x(t) \to 1 \,\,\, \mbox{as} \,\,\, t \to -\infty
$$
and
$$
x(t) \to 1 \,\,\, \mbox{or} \,\,\, x(t) \to -1 \,\,\, \mbox{as} \,\,\, t \to +\infty.
$$
More precisely, 
$$
x+1 \in H^{1}((-\infty,0]) \,\,\, \mbox{or} \,\,\, x-1 \in H^{1}((-\infty,0])
$$
and
$$
x+1 \in H^{1}([0,+\infty)) \,\,\, \mbox{or} \,\,\, x-1 \in H^{1}([0,+\infty)).
$$
\end{lem}

\noindent {\bf Proof.} \, Suppose the lemma is false. Then, there is $x \in H^{1}_{loc}(\mathbb{R})$ with $J_\epsilon(x)<\infty$, $\eta >0$ and a sequence $(t_n)$ with $|t_n| \to +\infty$ as $n \to +\infty$ such that
\begin{equation} \label{E3}
x(t_n) \in (-\infty, -1-\eta) \cup (-1+\eta, 1 - \eta) \cup (1+\eta, \infty).
\end{equation} 
Let
\begin{equation} \label{E4}
d=\{\tilde{V}(x): x \in (-\infty, -1- {\eta}/{2}) \cup (-1 + {\eta}/{2},1- {\eta}/{2}) \cup (1+{\eta}/{2}, +\infty)\}>0.
\end{equation}
We can assume, without loss of generality, $t_n \to +\infty$ and $t_{n+1} \geq t_n +1$ for all $n \in \mathbb{N}$. If 
$$
x(t) \in(-\infty, -1- {\eta}/{2}) \cup (-1 + {\eta}/{2},1- {\eta}/{2}) \cup (1+{\eta}/{2}, +\infty) \,\,\, \forall t \in [t_n,t_{n+1}],
$$
we have that
\begin{equation} \label{E5}
\int_{t_n}^{t_n +1}a(\epsilon t)\tilde{V}(x(t))\,dt \geq  l_0d \,\,\, \forall n \in \mathbb{N}.
\end{equation}
Otherwise, there exists $t^{*} \in [t_n, t_n +1]$ with $|x(t^*)-x(t_n)| \geq \eta / 2$. Thereby, 
$$
\eta / 2 \leq |x(t^*)-x(t_n)| \leq \int_{t_n}^{t^{*}}|\dot{x}|\,dt  \leq \sqrt{t^{*}-t_n}\left( \int_{t_n}^{t^{*}}|\dot{x}|^{2}\, dt \right)^{\frac{1}{2}},
$$
from where it follows that
\begin{equation} \label{E6}
 \int_{t_n}^{t_n+1}|\dot{x}|^{2}\, dt \geq \eta^{2} / 4.
\end{equation}
From (\ref{E5}) and (\ref{E6}), 
$$
 \int_{t_n}^{t_n+1}\left(\frac{1}{2}|\dot{x}|^{2}+a(\epsilon t)\tilde{V}(x(t)) \right)\,dt \geq \min\{l_0d, \eta^{2} / 4\}
$$
and so,
$$
J_\epsilon(x) \geq \sum_{n=1}^{+\infty} \int_{t_n}^{t_n+1}\left(\frac{1}{2}|\dot{x}|^{2}+a(\epsilon t)\tilde{V}(x(t)) \right)\,dt = +\infty
$$
which is  a contradiction, because by hypothesis $J_{\epsilon}(x)< \infty$. Then, 
$$
x(t) \to -1 \,\,\, \mbox{or} \,\,\, x(t) \to 1 \,\,\, \mbox{as} \,\,\, t \to +\infty.
$$ 
The same argument works to prove that  
$$
x(t) \to -1 \,\,\, \mbox{or} \,\,\, x(t) \to 1 \,\,\, \mbox{as} \,\,\, t \to -\infty. 
$$
By (\ref{EST1})-(\ref{EST2}), if $x(t) \to 1$ as $t \to +\infty$, there is $T>0$ such that
$$
\int_{T}^{+\infty}(x(t)-1)^{2} \leq \int_{T}^{+\infty} \frac{\tilde{V}(x(t))}{C_1}\,dt \leq \frac{1}{l_0 C_1} \int_{T}^{+\infty} a(\epsilon t)\tilde{V}(x(t))\,dt \leq \frac{1}{l_0 C_1} J_\epsilon(x) < \infty.
$$
The above inequality yields $x-1 \in H^{1}([0,+\infty))$. Analogous approach can be repeated to the cases 
$$
x(t) \to -1 \,\,\, \mbox{as} \,\,\ t \to +\infty, \,\, x(t) \to 1 \,\,\, \mbox{as} \,\,\ t \to -\infty \,\,\, \mbox{and}\,\,\, x(t) \to -1 \,\,\, \mbox{as} \,\,\ t \to -\infty.
$$
$\hfill \rule{2mm}{2mm}$

The next lemma will be used to study the convergence of the Palais-Smale sequences associated with $J_\epsilon$.  

\begin{lem} \label{L2} Let $A,T>0$. There is $B>0$, such that if $x \in H^{1}_{loc}(\mathbb{R})$ with $J_\epsilon(x) \leq A$, then $\|x\|_{H^{1}([-T,T])} \leq B$.

\end{lem}

\noindent {\bf Proof.} \, First of all, note that  
$$
\int_{-T}^{T}|\dot{x}|^{2}\,dt \leq 2A.
$$
By coercivity of $\tilde{V}$, there exists $C>0$ such that 
$$
\tilde{V}(x) > \frac{A}{l_0T} \,\,\,\mbox{for} \,\,  |x| \geq C. 
$$
Once
$$
\int_{-T}^{T}a(\epsilon t) \tilde{V}(x(t))\,dt \leq A,
$$
there is $t^{*} \in [-T,T]$ such that $\tilde{V}(x(t^{*})) \leq \frac{A}{2T}$ and $|x(t^{*})| \leq C$. Hence, for all $s \in [-T,T]$,
$$
|x(s)| \leq |x(t^{*})| + \left| \int_{\min\{s,t^*\}}^{\max\{s,t^{*}\}}\dot{x}(t)\, dt \right| \leq |x(t^{*})| + \sqrt{|s-t^{*}|}\left| \int_{\min\{s,t^*\}}^{\max\{s,t^{*}\}}|\dot{x}|^{2}\,dt \right|^{\frac{1}{2}} \leq C + 2\sqrt{TA},
$$ 
showing that
$$
\|x\|_{\infty} \leq C + 2\sqrt{AT}.
$$
$\hfill \rule{2mm}{2mm}$

\section{Existence of solution for Rabinowitz's condition}

In this section, we intend to prove Theorem \ref{T1}, by assuming that $a$ verifies the Rabinowitz's condition. 

\vspace{0.5 cm}

In what follows, we denote by ${\cal{B}_\epsilon}$, ${\cal{B}}_0$ and ${\cal{B}}_\infty$ the following real numbers
$$
{\cal{B}_\epsilon}=\inf \{J_\epsilon(x): \, x \in W \},
$$
$$
{\cal{B}}_0=\inf \{J_0(x): \, x \in W \},
$$
and
$$
{\cal{B}}_\infty=\inf \{J_\infty(x): \, x \in W \},
$$
where $ J_\infty:H^{1}_{loc}(\mathbb{R}) \to [0,+\infty]$ is the functional given by  
$$
J_\infty(x)=\int_{-\infty}^{+\infty}\left( \frac{1}{2}|\dot{x}|^{2}+a_\infty \tilde{V}(x(t))  \right)dt.
$$

Related to the above numbers, we have the ensuing result

\begin{lem} \label{031} The numbers ${\cal{B}_\epsilon}$, ${\cal{B}}_0$ and ${\cal{B}}_\infty$ verify
$$
{\cal{B}}_0 < {\cal{B}}_\infty \,\,\, \mbox{and} \,\,\, \lim_{\epsilon \to 0} {\cal{B}_\epsilon} = {\cal{B}}_0. 
$$
\end{lem}
\noindent {\bf Proof.} \,  In what follows, we denote by $w_0, w_\infty \in W$ the functions that verify
$$
\left\{
\begin{array}{l} 
\ddot{w_0}(t)=a(0) V'(w_0(t)), \,\,\, t \in \mathbb{R}, \\
w_0(t) \in (-1,1) \,\,\,\, \forall t \in \mathbb{R}, \\
w_0(t) \to -1 \,\,\, \mbox{as} \,\,\, t \to -\infty, \,\, w_0(t) \to 1 \,\,\, \mbox{as} \,\,\, t \to +\infty,
\end{array} 
\right.
\eqno{(P_0)}
$$
and
$$
\left\{
\begin{array}{l} 
\ddot{w_\infty}(t)=a_\infty V'(w_\infty(t)), \,\,\, t \in \mathbb{R}, \\
w_\infty(t) \in (-1,1) \,\,\,\, \forall t \in \mathbb{R}, \\ 
w_\infty(t) \to -1 \,\,\, \mbox{as} \,\,\, t \to -\infty, \,\, w_\infty(t) \to 1 \,\,\, \mbox{as} \,\,\, t \to +\infty.
\end{array}
\right.
\eqno{(P_\infty)}
$$
with 
$$
{\cal{B}}_0=J_0(w_0) \,\,\,  \mbox{and} \,\,\, {\cal{B}}_\infty=J_\infty(w_\infty).
$$ 
The existence of $w_0$ and $w_\infty$ was established in \cite{BS}.

By hypothesis $a_0 < a_\infty$, then
$$
{\cal{B}}_0 \leq J_0(w_\infty) < J_\infty(w_\infty)= {\cal{B}}_\infty,
$$
showing the first part of the lemma. For the second part, we begin observing  that
$$
J_0(w) \leq J_\epsilon(w) \,\,\, \forall  w \in W.
$$
Consequently, 
$$
{\cal{B}}_0 \leq {\cal{B}_\epsilon} \,\,\, \forall \epsilon >0,
$$
leading to 
\begin{equation} \label{E7}
{\cal{B}}_0 \leq \liminf_{\epsilon \to 0}{\cal{B}_\epsilon}.
\end{equation}
On the other hand, since $w_0 \in W$, 
$$
{\cal{B}_\epsilon} \leq J_\epsilon(w_0)=\int_{-\infty}^{+\infty}\left( \frac{1}{2}|\dot{w_0}|^{2}+a(\epsilon t)\tilde{V}(w_0(t))  \right)dt.
$$
Using Lebesgue's Theorem, we deduce that
\begin{equation} \label{E8}
\limsup_{\epsilon \to 0}{\cal{B}_\epsilon} \leq \int_{-\infty}^{+\infty}\left( \frac{1}{2}|\dot{w_0}|^{2}+a(0)\tilde{V}(w_0(t))  \right)dt=J_0(w_0)={\cal{B}}_0.
\end{equation}
From (\ref{E7})-(\ref{E8}), 
$$
\limsup_{\epsilon \to 0}{\cal{B}_\epsilon}= {\cal{B}}_0.
$$
$\hfill \rule{2mm}{2mm}$

The next lemma establishes that minimum points of $J$ on $W$ are in fact solutions for (\ref{E1})-(\ref{E2}).

\begin{lem} \label{03} If $x \in W$ verifies $J_\epsilon(x)={\cal{B}}_\epsilon$, then $x$ solves problem (\ref{E1})-(\ref{E2}) and \linebreak $x(t) \in (-1,1)$ for all $t \in \mathbb{R}$.

\end{lem}

\noindent {\bf Proof.} \, We start the proof recalling that 
$$
x+hv \in W \,\,\, \mbox{for all} \,\,\, v \in H^{1}(\mathbb{R}) \,\,\, \mbox{and} \,\,\, h \in \mathbb{R}.
$$
Since $J_\epsilon(x)={\cal{B}}_\epsilon$, the above information yields
$$
\frac{J_\epsilon(x+hv)-J_\epsilon(x)}{h} \geq 0 \,\,\, \forall h \geq 0.
$$
Letting the limit of $h \to 0$, we get
$$
J_\epsilon'(x)v \geq 0 \,\,\, \forall v \in H^{1}(\mathbb{R})
$$
and so,
$$
J_\epsilon'(x)v =0  \,\,\, \forall v \in H^{1}(\mathbb{R}), 
$$
implying that $x$ is a critical point of $J_\epsilon$. Therefore, $x$ is a solution of O.D.E.
$$
\ddot{x}(t)=a(\epsilon t)\tilde{V}'(x(t)), \,\,\,\, t \in \mathbb{R}.
$$
Moreover, by $x \in W$, one have 
$$
x(t) \to -1 \,\,\, \mbox{as} \,\,\, t \to -\infty, \,\, x(t) \to 1 \,\,\, \mbox{as} \,\,\, t \to +\infty.
$$ 

Now, we will prove that 
$$
x(t) \in (-1,1) \,\, \forall \,\, t \in \mathbb{R}.
$$   
If $x(t) >1$ for some $t \in \mathbb{R}$, then let $t_0 \in \mathbb{R}$ with $x(t_0)=\displaystyle \max_{t \in \mathbb{R}}x(t)>1$. Thereby, 
$$
\ddot{x}(t_0) \leq 0 \,\,\, \mbox{and} \,\,\, \tilde{V}'(x(t_0)) >0,
$$
which is an absurd. Thus $x(t) < 1$ for all $t \in \mathbb{R}$. The same type of argument works to show that $x(t) <-1$ for all $t \in \mathbb{R}$. From the above information, we can conclude that $x$ is a solution for original problem (\ref{E1})-(\ref{E2}), because
$$
\tilde{V}(x(t))=V(x(t)) \,\,\, \forall t \in \mathbb{R},
$$
finishing the proof of lemma. $\hfill \rule{2mm}{2mm}$

\vspace{0.5 cm}

The next result shows that associated with ${\cal{B}_\epsilon}$, we have a Palais-Smale sequence for $J$.

\begin{lem} \label{L031} 
There is a $(PS)_{\cal{B}_\epsilon}$ sequence for $J_\epsilon$. 
\end{lem}

\noindent {\bf Proof.}  Since $J_\epsilon$ is bounded from below, there is $(x_n) \subset W$ such that
$$
J_\epsilon(x_n) \to {\cal{B}_\epsilon} \,\,\,\, \mbox{as} \,\,\,\, n \to +\infty.
$$   
Now, it is easy to check that if $x,z \in W$, then $x-z \in H^{1}(\mathbb{R})$. Therefore, we can define on $W$ the metric $\rho:W \times W \to [0+\infty)$ given by 
$$
\rho(x,z)=\|x-z\|,
$$
where $\|\,\,\ \|$ denotes the usual norm in $H^{1}(\mathbb{R})$. A direct computation gives that $(W,\rho)$ is a complete metric space . Once $J_\epsilon$ is lower semicontinuous and bounded from below on $(W, \rho)$, by Ekeland's Variational Principle there is $(u_n) \subset W$ verifying
$$
\|x_n-u_n\|=o_{n}(1),
$$
$$
J_\epsilon(u_n) \to {\cal{B}_\epsilon},
$$
and
$$
J_\epsilon(u_n) \leq J_\epsilon(x)+\frac{1}{n}\rho(u_n,x) \,\,\, \forall x \in W,
$$
that is,
$$
J_\epsilon(u_n) \leq J_\epsilon(x) + \frac{1}{n}\|u_n-x\| \,\,\, \forall x \in W.
$$
Now, for each $v \in H^{1}(\mathbb{R})$ and $t \in (0,+\infty)$, we know that
$$
u_n +tv \in W,
$$
then
$$
\frac{J_\epsilon(u_n+tv)-J_\epsilon(u_n)}{t} \geq -\frac{1}{n}\|v\|.
$$
Thus, taking the limit of $n \to +\infty$,  we get
$$
J'_\epsilon(u_n)v \geq -\frac{1}{n}\|v\|.
$$
From this,
$$
\|J_\epsilon(u_n)\|_* \leq \frac{1}{n},
$$ 
from where it follows that 
$$
\|J_\epsilon(u_n)\|_* \to 0 \,\,\, \mbox{as} \,\,\, n \to +\infty,
$$
showing that $(u_n)$ is a $(PS)_{\cal{B}_\epsilon}$ sequence for $J_\epsilon$. $\hfill \rule{2mm}{2mm}$

\vspace{0.5 cm}

The next lemma is crucial in our approach and its proof can be found in \cite{S}. 

\begin{lem} \label{L04} Let $x_0, x_1 \in (-1,1), x_0<x_1, t_0<t_1$ and $x \in H^{1}([t_0,t_1])$ with $x(t_0)=x_0$ and $x(t_1)=x_1$. Then, 
$$
\int_{t_0}^{t_1}\left( \frac{1}{2}|\dot{x}|^{2}+a_\infty V(x(t))\right)\,dt \geq \int_{w_\infty^{-1}(x_0)}^{w_\infty^{-1}(x_1)}\left( \frac{1}{2}|\dot{w_\infty}|^{2}+a_\infty V(w_\infty(t))\right)\,dt, 
$$
where $w_\infty $ was given in the proof of Lemma \ref{031}. 
\end{lem}

\vspace{0.5 cm}

The main result this section can be stated as follows 

\begin{thm} \label{T12} Assume that $(V_1)-(V_4)$ hold. If $a$ belongs to Class 8, there is $\epsilon^* >0$, such that problem (\ref{E1})-(\ref{E2}) has a solution $x \in H_{loc}^{1}(\mathbb{R}) \cap C^{2}(\mathbb{R})$ for all $\epsilon \in (0, \epsilon^*)$. Moreover, $x(t) \in (0,1)$ for all $t \in \mathbb{R}$.

\end{thm}

\noindent {\bf Proof.} \, First of all, by Lemma \ref{03}, we see that to prove Theorem \ref{T1}, it is enough to show that there exists $\epsilon^*>0$ such that ${\cal{B}}_\epsilon$ is achieved for all $\epsilon \in [0, \epsilon^*)$. 

To prove that ${\cal{B}}_\epsilon$ is achieved, we begin  recalling that from Lemma \ref{L031}, there is a $(PS)_{{\cal{B}}_\epsilon}$ sequence for $J_\epsilon$, that is, there exists $(x_n) \subset W$ such that
$$
J_\epsilon(x_n) \to {\cal{B}}_\epsilon \,\,\, \mbox{and} \,\,\, \|J_\epsilon'(x_n)\|_* \to 0 \,\,\, \mbox{as} \,\,\, n \to +\infty.
$$ 
From this, 
$$
J_\epsilon(x_n) \leq A= \sup_{n}J_\epsilon(x_n) \,\,\, \forall n \in \mathbb{N}.
$$
By Lemma \ref{L2}, for each $T>0$, there is $B=B(T,A)>0$ such that 
$$
\|x_n\|_{H^{1}([-T,T])} \leq B \,\,\ \forall n \in \mathbb{N}.
$$
Hence, there is a subsequence of $(x_n)$, still denoted by itself, and $x \in H^{1}_{loc}(\mathbb{R})$ verifying 
$$
x_n \to x \,\,\, \mbox{uniformly in} \,\, [-T,T] \,\,\, \mbox{and} \,\,\, x_n \rightharpoonup x \,\,\, \mbox{in} \,\,\, H^{1}([-T,T])\,\,\, \forall T >0. 
$$
Combining these limits with the fact that $J_\epsilon$ is lower semicontinuous, we also derive that
\begin{equation} \label{Z01}
J_\epsilon(x) \leq  {\cal{B}}_\epsilon.
\end{equation}
Next, we will show that $J'_\epsilon(x)=0$. To see why, note that for each $v \in C_{0}^{\infty}(\mathbb{R})$ fixed, we have that $J_\epsilon'(x_n)v=o_n(1)$. Then, 
$$
\int_{\alpha}^{\beta}\dot{x_n}\dot{v}\, dt + \int_{\alpha}^{\beta}a(\epsilon t)\tilde{V}'(x_n(t))v(t)\,dt=o_n(1),
$$
where $supp \, v \subset [\alpha, \beta]$. Letting $n \to +\infty$, we get
$$
\int_{\alpha}^{\beta}\dot{x}\dot{v}\, dt + \int_{\alpha}^{\beta}a(\epsilon t)\tilde{V}'(x(t))v(t)\,dt=0,
$$
implying that $x$ is a solution of equation O.D.E.
$$
\ddot{x}(t)=a(\epsilon t)\tilde{V}'(x(t)),
$$
and so,
$$
J'_\epsilon(x)=0.
$$
Moreover, by Fatous' Lemma $J_\epsilon(x)<+\infty $. Consequently, by Lemma \ref{L1}
$$
x(t) \to -1 \,\,\, \mbox{or} \,\,\, x(t) \to 1 \,\,\, \mbox{as} \,\,\, t \to -\infty
$$
or
$$
x(t) \to 1 \,\,\, \mbox{or} \,\,\, x(t) \to -1 \,\,\, \mbox{as} \,\,\, t \to +\infty.
$$
Our next step is showing that below limit 
\begin{equation} \label{E18}
x(t) \to -1 \,\,\, \mbox{as} \,\,\, t \to +\infty
\end{equation}
does not hold. To this end, we suppose by contradiction that it holds and we will set for each $\tau>0$ the real number
$$
{\Lambda}_{\tau}=\int_{w_\infty^{-1}(-1+\tau)}^{w_\infty^{-1}(1-\tau)}\left( \frac{1}{2}|\dot{w_\infty}|^{2}+a_\infty V(w_\infty(t))\right)\,dt,
$$ 
where $w_\infty \in W$ and $J_\infty(w_\infty)={\cal{B}}_\infty$. By a routine calculus, 
\begin{equation} \label{E188}
{\Lambda}_{\tau} \to {\cal{B}}_\infty \,\,\, \mbox{as} \,\,\, \tau \to 0.
\end{equation}
In the last limit, we have used that $\displaystyle \lim_{t \to +\infty}w_\infty(t)=1$ and $\displaystyle \lim_{t \to -\infty}w_\infty(t)=-1$. 

The inequality ${\cal{B}}_0 < {\cal{B}}_\infty$ in conjunction with (\ref{E188}) implies that there is $\tau >0 $ small enough verifying
$$
\left( \frac{a_\infty- \tau}{a_\infty} \right){\Lambda}_{\tau} >{\cal{B}}_0.
$$
Now, by $(a_7)$, let $T>0$ be large enough so that $a(\epsilon t) \geq a_\infty- \tau$ on $[T,+\infty)$ and $x(T) < -1+\tau$. Let $n$ be a large enough that $x_n(T)<-1+\tau$. Let $T < \alpha < \beta$ with $x_n(\alpha)=-1+\tau$ and $x_n(\beta)=1-\tau$. By Lemma \ref{L04},
$$
J_\epsilon(x_n) \geq \left(\frac{a_\infty-\tau}{a_\infty}\right)\int_{\alpha}^{\beta}\left(\frac{1}{2}|\dot{x}_n|^{2}+a_\infty V(x_n(t)) \right)\,dt \geq \left(\frac{a_\infty-\tau}{a_\infty}\right){\Lambda}_{\tau}
$$
and so,
$$
{\cal B}_\epsilon = \lim_{n \to +\infty}J_\epsilon(x_n) \geq \left(\frac{a_\infty-\tau}{a_\infty}\right){\Lambda}_{\tau} > {\cal B}_0.
$$
Consequently,
$$
\lim_{\epsilon \to 0}{\cal B}_\epsilon \geq \left(\frac{a_\infty-\tau}{a_\infty}\right){\Lambda}_{\tau} > {\cal B}_0,
$$
contradicting Lemma \ref{031}. This way,
$$
x(t) \to 1 \,\,\, \mbox{as} \,\,\, t \to +\infty.
$$
A similar argument can be used to show that
$$
x(t) \to -1 \,\,\, \mbox{as} \,\,\, t \to -\infty.
$$
As in the proof of Lemma \ref{L1}, we derive that $x+1 \in H^{1}((-\infty,0])$ and $x-1 \in H^{1}([0,+\infty)$. Then,  $x \in W$, and by (\ref{Z01})$, J_\epsilon(x)={\cal B}_\epsilon$ finishing the proof. $\hfill \rule{2mm}{2mm}$

\section{Existence of solution for the asymptotically periodic case}

In this section, we intend to prove the existence of solution for (\ref{E1})-(\ref{E2}), by assuming that $a$ is asymptotically periodic. 

\vspace{0.5 cm}

The main result in section is the following

\begin{thm} \label{T2} Assume that $(V_1)-(V_4)$ hold. If $a$ belongs to Class 9,  problem (\ref{E1})-(\ref{E2}) has a solution $x \in H_{loc}^{1}(\mathbb{R}) \cap C^{2}(\mathbb{R})$ for each $\epsilon >0$. Moreover, $x(t) \in (0,1)$ for all $t \in \mathbb{R}$.

\end{thm}

In the proof of Theorem \ref{T2}, without loss of generality, we assume that $\epsilon=1$. Moreover,  we will use the fact that problem (\ref{E1})-(\ref{E2}) has an increasing solution $w_P \in H_{loc}^{1}(\mathbb{R}) \cap C^{2}(\mathbb{R})$ with $w_P \in W$ and $J_P(w_P)={\cal {B}}_P$, where  $ J_P:H^{1}_{loc}(\mathbb{R}) \to [0,+\infty]$ is the functional given by  
$$
J_P(x)=\int_{-\infty}^{+\infty}\left( \frac{1}{2}|\dot{x}|^{2}+a_P(t) \tilde{V}(x(t))  \right)dt.
$$
and 
$$
{\cal {B}}_P= \inf \{J_P(x): \, x \in W \}.
$$
The existence of $w_P$ can be seen in \cite{BS}.

In the sequel, we denote by $ J:H^{1}_{loc}(\mathbb{R}) \to [0,+\infty]$ the functional given by  
$$
J(x)=\int_{-\infty}^{+\infty}\left( \frac{1}{2}|\dot{x}|^{2}+a(t) \tilde{V}(x(t))  \right)dt
$$
and by ${\cal {B}}$, the real number given by
$$
{\cal {B}}= \inf \{J(x): \, x \in W \}.
$$
Here, we would like point out that all results proved in Section 2 are true for functionals $J$ and $J_P$.  Moreover, from $(a_9)$, we also have 
\begin{equation} \label{E10}
{\cal {B}} < {\cal {B}}_P.
\end{equation}

\vspace{0.5 cm}

\noindent {\bf Proof of Theorem \ref{T2}}

\vspace{0.5 cm}
As in the proof of Theorem \ref{T12}, our main goal is to show that $\cal{B}$ is achieved on $W$. 
Hereafter, $a_0=\displaystyle \inf_{t \in \mathbb{R}}a(t)$ and we fix $\delta>0$ such that
\begin{equation} \label{E11}
{\cal {B}} + \delta < {\cal {B}}_P. 
\end{equation}
Moreover, we also fix $M=M(\delta)>0$ such that
\begin{equation} \label{E12}
|a(t)-a_P(t)|< \frac{\delta a_0}{2{\cal {B}}} \,\,\, \mbox{for} \,\,\, |t|>M
\end{equation}
and $\epsilon >0$ verifying
\begin{equation} \label{E13}
V(z) < \frac{\delta}{4M\|a_P\|_\infty} \,\,\, \forall z \in [-1,-1+\epsilon/2] \cup [1-\epsilon/2, 1].
\end{equation}
For $\epsilon >0$ given above, combining the same arguments explored in \cite{BS}  with Ekeland's variational principle, we can find sequences $(U_n) \subset W$, $(s_n),(t_n) \subset \mathbb{R}$ with $s_n < t_n$ satisfying:
$$
J(U_n) \to {\cal {B}}, \,\, J'(U_n) \to 0  \,\,\, \mbox{as} \,\,\, n \to +\infty,
$$
$$
U_n(t) \in [-1,-1+\epsilon/2] \,\,\, \forall t \in (-\infty, s_n], 
$$ 
$$
U_n(t) \in [1-\epsilon/2,1] \,\,\, \forall t \in [t_n, +\infty),
$$
$$
U_n(t) \in [-1+\epsilon/2 , 1-\epsilon/2] \,\,\, \forall t \in [s_n,t_n],
$$
$$
U_n(t_n)=1-\epsilon +o_{n}(1), \,\,\, U_n(s_n)=-1+\epsilon +o_{n}(1),
$$
and 
$$
(t_n-s_n) \,\,\,\, \mbox{is bounded in } \,\,\, \mathbb{R}.
$$
A direct computation shows that for some subsequence of $(U_n)$, still denoted by itself, there is $U \in C(\mathbb{R}) \cap H^{1}_{loc}(\mathbb{R})$ such that
\begin{equation} \label{E131}
U_n \to U \,\,\, \mbox{in} \,\,\, C_{loc}(\mathbb{R}).
\end{equation}
As in the proof Theorem \ref{T12}, we see that
$$
J(U) \leq {\cal{B}} \,\,\,\, \mbox{and} \,\,\,\, J'(U)=0. 
$$
This way, the theorem follows provided that $U \in W$. To show this fact, we make the following claim

\begin{claim} \label{C1}
The sequence $(s_n)$ is bounded. 
\end{claim}

\vspace{0.5 cm}

Indeed, if the claim is not true, we must have for some subsequence 
$$
s_n \to +\infty \,\,\,\, \mbox{or} \,\,\,\, s_n \to -\infty.
$$
Using the above limits, we deduce that
$$
U(t) \in [-1,-1+\epsilon/2] \cup [1-\epsilon/2,1] \,\,\, \forall t \in \mathbb{R}.
$$
Thus by (\ref{E13}),
\begin{equation} \label{E14}
\tilde{V}(U(t)) < \frac{\delta}{8M\|a_P\|_\infty} \,\,\, \forall t \in \mathbb{R}.
\end{equation}
Note that
\begin{equation} \label{E15}
J(U_n)=J_P(U_n)+\int_{\mathbb{R}}(a(t)-a_P(t))\tilde{V}(U_n(t))\,dt \geq {\cal {B}}_P + \int_{\mathbb{R}}(a(t)-a_P(t))\tilde{V}(U_n(t))\,dt.
\end{equation}
Since 
$$
\int_{\mathbb{R}}|a(t)-a_P(t)|\tilde{V}(U_n(t))\,dt = \int_{|t|\leq M}|a(t)-a_P(t)|\tilde{V}(U_n(t))\,dt + \int_{|t|>M}|a(t)-a_P(t)|\tilde{V}(U_n(t))\,dt,
$$
by (\ref{E12}) and (\ref{E131}),
$$
\limsup_{n \to +\infty}\int_{\mathbb{R}}|a(t)-a_P(t)|\tilde{V}(U_n(t))\,dt \leq \int_{|t|\leq M}|a(t)-a_P(t)|\tilde{V}(U(t))\,dt + \frac{\delta}{2}.
$$
Now, using (\ref{E14}), 
$$
\limsup_{n \to +\infty}\int_{\mathbb{R}}|a(t)-a_P(t)|\tilde{V}(U_n(t))\,dt < \frac{\delta}{2}+\frac{\delta}{2}=\delta.
$$
Thereby, there is $n_0 \in \mathbb{N}$ such that 
\begin{equation} \label{E16}
\int_{\mathbb{R}}|a(t)-a_P(t)|\tilde{V}(U_n(t))\,dt< \delta \,\,\, \forall n \geq n_0.
\end{equation}
Combining (\ref{E15}) with (\ref{E16}), we derive
$$
J(U_n)=J_P(U_n)+\int_{\mathbb{R}}(a(t)-a_P(t))\tilde{V}(U_n(t))\,dt \geq {\cal {B}}_P - \delta \,\,\, \forall n \geq n_0,
$$
that is, 
$$
J(U_n) \geq {\cal {B}}_P - \delta \,\,\, \forall n \geq n_0.
$$
Taking the limit of $n \to +\infty$ in the last inequality, we obtain the estimate
$$
{\cal {B}} > {\cal {B}}_P - \delta,
$$
which contradicts (\ref{E11}).

The boundedness of $(s_n)$ implies that $(t_n)$ is also bounded, thus we can assume without of generality, that there are $t,s \in \mathbb{R}$ verifying
$$
s_n \to s \,\,\,\, \mbox{and} \,\,\, t_n \to t \,\,\, \mbox{as} \,\,\, n \to +\infty.
$$ 
This way,
$$
U(z) \in [-1,-1+\epsilon] \,\,\, \forall z \in (-\infty,s] \,\,\, \mbox{and} \,\,\, U(z) \in [1-\epsilon,1] \,\,\, \forall z \in [t,+\infty).
$$
The above information together with Lemma \ref{L1} gives
$$
U(z) \to 1 \,\,\ \mbox{as} \,\,\, z \to +\infty \,\,\, \mbox{and} \,\,\,  U(z) \to -1 \,\,\ \mbox{as} \,\,\, z \to -\infty.
$$
This finishes the proof of Theorem \ref{T2}.  $\hfill \rule{2mm}{2mm}$

\section{Existence of solution for the coercive case}

In this section, we intend to prove the existence of solution for (\ref{E1})-(\ref{E2}), by assuming that $a$ is {\it coercive}. Here, our main result has the following statement

\begin{thm} \label{T3} Assume that $(V_1)-(V_4)$ hold. If $a$ is coercive, problem (\ref{E1})-(\ref{E2}) has a solution $x \in H_{loc}^{1}(\mathbb{R}) \cap C^{2}(\mathbb{R})$ for all  $\epsilon >0$. Moreover, $x(t) \in (0,1)$ for all $t \in \mathbb{R}$.

\end{thm}

In the sequel, we will assume that $\epsilon =1$. However, in the proof of the above result, we must to be careful to use the arguments of the previous sections, more precisely Section 2. In the sequel, we need to fix the following sets 
$$
W_a=\left\{ x \in H^{1}_{loc}(\mathbb{R}): x+1 \in H_a^{1}((-\infty, 0]), \, x-1 \in H_a^{1}([0,+\infty))  \right\}
$$
where
$$
H_a^{1}((-\infty, 0])=\left\{v \in H^{1}((-\infty, 0])\,:\, \int_{-\infty}^{0}a(t)|v(t)|^{2}\,dt<+\infty \right\}
$$
endowed with the norm
$$
\|v\|_{a,-\infty}=\left( \int_{-\infty}^{0}|v'(t)|^{2}\,dt+\int_{-\infty}^{0}a(t)|v(t)|^{2}\,dt \right)^{\frac{1}{2}}.
$$
The space $ H_a^{1}([0, +\infty))$ is defined of a similar way, that is,
$$
H_a^{1}([0, +\infty))=\left\{v \in H^{1}([0, +\infty))\,:\, \int_{0}^{+\infty}a(t)|v(t)|^{2}\,dt<+\infty \right\},
$$
endowed with the norm
$$
\|v\|_{a,+\infty}=\left( \int_{0}^{+\infty}|v'(t)|^{2}\,dt+\int_{0}^{+\infty}a(t)|v(t)|^{2}\,dt \right)^{\frac{1}{2}}.
$$
From $(a_{10})$, we know that $\displaystyle \inf_{t\in \mathbb{R}}a(t)>0$, then the below embeddings are continuous
$$
H_a^{1}([0, +\infty)) \hookrightarrow H^{1}([0, +\infty)) \,\,\, \mbox{and} \,\,\, H_a^{1}((-\infty, 0]) \hookrightarrow H^{1}((-\infty, 0]).
$$

\vspace{0.5 cm}

\noindent {\bf Proof Theorem \ref{T3}}

\vspace{0.5 cm}

Hereafter, we follow the same approach of the previous section. Adapting the same arguments explored in \cite{BS}, we can find sequences $(U_n) \subset W_a$, 
$(s_n),(t_n) \subset \mathbb{R}$ with $s_n < t_n$ satisfying:
$$
J(U_n) \to {\cal {B}}, \, J'(U_n) \to 0 \,\,\,\, \mbox{as} \,\,\, n \to +\infty,
$$
$$
U_n(t) \in [-1,-1+\epsilon/2] \,\,\, \forall t \in (-\infty, s_n], 
$$ 
$$
U_n(t) \in [1-\epsilon/2,1] \,\,\, \forall t \in [t_n, +\infty),
$$
$$
U_n(t) \in [-1+\epsilon/2, 1-\epsilon/2] \,\,\, \forall t \in [s_n,t_n],
$$
$$
U_n(t_n)=1-\epsilon/2+o_{n}(1), \,\,\, U_n(s_n)=-1+\epsilon/2+o_n(1),
$$
and 
$$
(t_n-s_n) \,\,\,\, \mbox{is bounded in } \,\,\, \mathbb{R}.
$$
A direct computation shows that for some subsequence of $(U_n)$, still denoted by itself, there is $U \in C(\mathbb{R})\cap H_{loc}^{1}(\mathbb{R})$ such that
\begin{equation} \label{E1311}
U_n \to U \,\,\, \mbox{in} \,\,\, C_{loc}(\mathbb{R}).
\end{equation}
Moreover,
$$
J(U) \leq {\cal {B}} \,\,\, \mbox{and} \,\,\,\, J'(U)=0.
$$
Here, $J$ and ${\cal {B}}$ are as in the proof of Theorem \ref{T2}. Our goal is proving that ${\cal {B}}$ is achieved on $W$. To this end, we will study again the behavior of sequence $(s_n)$.

\vspace{0.5 cm}

\begin{claim} \label{C2}
The sequence $(s_n)$ is bounded. 
\end{claim}

\vspace{0.5 cm}

Arguing by contradiction, we will assume that $(s_n)$ is unbounded. Then for some subsequence, still denoted by itself, we have that
\begin{equation} \label{E17}
s_n \to +\infty \,\,\, \mbox{or} \,\,\,\, s_n \to -\infty.
\end{equation}
Using the definition of $J$ and the properties of $(U_n)$, we derive that
$$
J(U_n) \geq V_0A_n(t_n-s_n)
$$
where
$$
A_n=\min_{z \in [s_n,t_n]}a(z) \,\,\,\, \mbox{and} \,\,\,\, V_0=\min_{-1+\epsilon/2 \leq z \leq 1-\epsilon/2}V(z)>0.
$$
Using the fact that $a$ is coercive in conjunction with (\ref{E17}) and the boundedness of $((t_n-s_n))$, we deduce that
$$
A_n \to +\infty \,\,\, \mbox{as} \,\,\, n \to +\infty.
$$
Since $(J(U_n))$ is bounded, the last inequality implies that
\begin{equation} \label{E18}
t_n-s_n \to 0 \,\,\, \mbox{as} \,\,\, n \to +\infty.
\end{equation}
Once $(U_n) \subset H_{loc}^{1}(\mathbb{R})$, for all $s,t \in \mathbb{R}$ the below inequality occurs 
$$
|U_n(t)-U_n(s)|\leq \sqrt{|t-s|}\left( \int_{\min\{t,s\}}^{\max\{t,s\}}|U_n'(r)|^{2}dr \right)^{\frac{1}{2}} \,\,\,\,\, \forall n \in \mathbb{N}.
$$
Thus,
$$
|U_n(t)-U_n(s)|\leq \sqrt{2|t-s|}J(U_n)^{\frac{1}{2}} \,\,\,\,\, \forall n \in \mathbb{N}.
$$
Now, the  boundedness of $(J(U_n))$ together with (\ref{E18}) gives
$$
|U_n(t_n)-U_n(s_n)| \to 0 \,\,\,\, n \to +\infty.
$$
However, this limit cannot occur, because
$$
|U_n(t_n)-U_n(s_n)|=2-\epsilon + o_n(1) \,\,\, \forall n \in \mathbb{N}.
$$
Therefore, the Claim \ref{C2} is proved. 

Now, the proof of Theorem \ref{T3} follows the same steps of the proof of Theorem \ref{T2}.

$\hfill \rule{2mm}{2mm}$

\section{Final remarks} In Section 2, we can remove the condition that $a \in L^{\infty}(\mathbb{R})$. Hovewer, we must work with the same spaces used in Section 5.

\end{document}